\documentclass[]{amsart}

\usepackage{color}
\usepackage{graphicx,epstopdf}
\usepackage{amsmath}
\usepackage{amssymb}
\usepackage{amsfonts}
\usepackage{amsthm}

\theoremstyle{remark}
\newtheorem{remark}{Remark}[section]
 
\begin{document}
 \title[]{Spectral approach to Korteweg-de Vries equations on the 
 compactified real line}

\author{Christian Klein}
\address{Institut de Math\'ematiques de Bourgogne, UMR 5584\\
                Universit\'e de Bourgogne-Franche-Comt\'e, 9 avenue Alain Savary, 21078 Dijon
                Cedex, France\\
    E-mail Christian.Klein@u-bourgogne.fr}

\author{Nikola Stoilov}
\address{Institut de Math\'ematiques de Bourgogne, UMR 5584\\
                Universit\'e de Bourgogne-Franche-Comt\'e, 9 avenue Alain Savary, 21078 Dijon
                Cedex, France\\
    E-mail Nikola.Stoilov@u-bourgogne.fr}
\date{\today}

\begin{abstract}
  We present a numerical approach for generalised Korteweg-de Vries (KdV) equations on the real line. In the spatial dimension we  compactify the real line and apply a Chebyshev collocation method. The time integration is  performed with an implicit Runge-Kutta method of fourth order. Several examples are discussed: initial data bounded but not vanishing at infinity as well as data not satisfying the Faddeev condition, i.e. with a slow decay towards infinity. 
\end{abstract}

 
\thanks{This work is partially supported by 
the ANR-FWF project ANuI - ANR-17-CE40-0035, the isite BFC project 
NAANoD, the ANR-17-EURE-0002 EIPHI and by the 
European Union Horizon 2020 research and innovation program under the 
Marie Sklodowska-Curie RISE 2017 grant agreement no. 778010 IPaDEGAN}
\maketitle

\section{Introduction}

Generalised Korteweg-de Vries (gKdV) equations,
\begin{equation}
    u_{t}(x,t)+u_{xxx}(x,t)+u(x,t)^{p-1}u_{x}(x,t)=0,
    \label{gKdV}
\end{equation}
where $p\in \mathbb{N}$, $u:\mathbb{R}\times \mathbb{R}^{+}\mapsto 
\mathbb{R}$, appear as asymptotic models in hydrodynamics, nonlinear 
optics, plasma physics, Bose-Einstein condensates, and   
essentially in most situations where predominantly 
one-dimensional phenomena are discussed and where dispersion 
dominates dissipation. Whereas this applies in particular to the case 
of the classical Korteweg-de Vries (KdV) equation ($p=2$), there are applications for instance 
in electrodynamics for the modified KdV equation ($p=3$), see 
\cite{mkdv}. Because of 
their importance in applications, there has been a considerable 
activity in developing numerical approaches for the gKdV equations. 
For initial data which are periodic or rapidly decreasing, numerical 
approaches based on the approximation of $u$ in (\ref{gKdV}) via 
truncated Fourier series, i.e., 
trigonometric polynomials, are very efficient, see for instance 
\cite{KLE08,KP13} and references therein.  


The Fourier approach, that is restricting the data to an interval 
$L[-\pi,\pi]$ ($L=const$, $L\in\mathbb{R}^{+}$) and continuing them periodically with 
period $2\pi L$ on the whole real line works very well for periodic 
functions and exhibits \emph{spectral convergence}, namely an 
exponential decrease of the numerical error with the number of 
Fourier modes. Schwartz class functions can be treated as periodic as 
one works with a finite precision and since $L$ can be chosen large 
enough such that all necessary derivatives of $u$ vanish at the domain 
boundaries with numerical precision. However, the same approach for 
initial data which do not tend to zero or are only slowly decreasing to zero for 
$x\to\infty$, would in general imply  a Gibbs phenomenon at the 
domain boundaries. The resulting method would therefore be only of 
first order in the number of Fourier modes, making it impossible to 
reach the high resolution necessary to treat e.g. rapid oscillations 
in the solution. Therefore, such situations are typically treated on a finite interval. This leads to the problem of how to impose boundary conditions, so that inside the computational domain the solution is the same as if the computation was done on the whole real line. 
Such boundary conditions are called \emph{transparent}.
B\'erenger \cite{Berenger} introduced \emph{perfectly 
matched layers} (PML) in electrodynamics to address this problem by 
extending the computational domain to layers glued to the domain 
boundaries. Inside the layers, the equation under consideration is deformed 
to a dissipative one, which is chosen to efficiently dissipate the 
solution. Whereas this works 
well for linear equations, examples for the nonlinear Schr\"odinger 
equation, see \cite{zhengpml,birem}, showed that in a nonlinear 
setting there will be back
reflections from the layers to the computational domain. For integrable equations exact transparent boundary 
conditions (TBC) can be given, e.g. for the case of modified KdV see \cite{zhengmkdv} based 
on \cite{fokasmkdv}. The problem with both PML and TBC is that they 
in general require initial data with compact support within the 
computational domain, thus limiting the class of solutions that can be studied.
The goal of the present paper is to establish a spectral numerical 
approach for generalised KdV 
equations (\ref{gKdV}) for initial data that are analytic on the whole real line, 
that are slowly decreasing towards infinity, or are bounded there. This numerical approach exhibits spectral convergence on the whole 
real line with a technique similar to \cite{orszag}.

Both the classical KdV equation and the modified KdV equation are completely integrable, which means 
 they have an infinite number of conserved quantities (for a comprehensive review on integrability see \cite{ZAK91}). In all 
other cases, the generalised KdV equations have only three conserved 
quantities, $\int_{\mathbb{R}}^{}u dx$ and the $L^{2}$ norm of $u$ and the energy
\begin{equation}
    E[u] = \int_{\mathbb{R}}^{} \left( 
    \frac{u^{p+1}}{p(p+1)}-\frac{1}{2} u_{x}^{2}\right)dx
    \label{energy}.
\end{equation}
The complete integrability of the classical KdV equation made it one of the best studied non-linear dispersive equations with a 
rather complete understanding of its solutions. However, even in this case 
there are open questions which motivate us to provide numerical tools 
to complement analytical studies in this context. The standard 
inverse scattering approach to KdV is only applicable if the 
Faddeev decay condition \cite{faddeev},
\begin{equation}
    \int_{\mathbb{R}}^{}(1+|x|)|u_{0}(x)|dx<\infty
    \label{faddeev},
\end{equation}
holds for the initial data $u(x,0)=u_{0}(x)$. 
  The direct scattering approach  to KdV  involves the determination of the spectrum of the 
Schr\"odinger equation for the potential $u_{0}(x)$,
$$\psi_{xx}+u_{0}(x)\psi=E\psi,$$ see \cite{ZAK91}.
It is known that this discrete spectrum is finite if the Faddeev condition 
(\ref{faddeev}) is satisfied, see \cite{Mar}. But except for the 
periodic case, see for instance \cite{algebro}, there is no complete 
understanding of the case of initial data not satisfying 
(\ref{faddeev}). The goal of the present paper is to 
provide numerical tools to study such cases. Of course, numerically one 
will be only able to study finite times and thus will not be able to 
address the question whether the time evolution of such data can 
lead to an infinite number of solitons.


The paper is organized as follows: In section \ref{Sec_Theory} we summarize a few 
theoretical facts on generalised KdV equations. In section \ref{Sec_Num_appr} we choose a 
compactification of the real line and describe the numerical approach for 
the generalised KdV equations. In section \ref{Sec_Examples} we discuss several examples. 
Concluding remarks are added in section \ref{Sec_Conclusion}.

\section{Theoretical preliminaries}\label{Sec_Theory}
In this section we summarize basic facts about generalised KdV 
equations needed in the following.

Though the KdV equations (\ref{gKdV}) are only completely integrable 
for $p=2$ and $p=3$, they have for all integer values of $p\geq2$ a 
solitary travelling wave solution which is explicitly given by 
$u=Q_{c}(x-x_{0}-ct)$ with $x_{0},c=const$ and with
\begin{equation}
    Q_{c}(z) = \left( 
    \frac{p(p+1)c}{2}\,\mbox{sech}^{2}\frac{\sqrt{c}(p-1)}{2}z\right)^{1/(p-1)}
    \label{soliton}.
\end{equation}
Thus, we have $Q_{c}(z)=c^{1/(p-1)}Q(\sqrt{c}z)$, where we have put 
$Q:=Q_{1}$. This simple scaling property of the \emph{solitons} allows to 
concentrate on the case $c=1$. If we refer in the following to 
the generalised KdV soliton, it is always implied that $c=1$. It was shown in \cite{BSW1987} that these solitons are linearly 
unstable for $p> 4$. Note that the energy of the soliton vanishes for 
$p=5$.

The generalised KdV equation has the following scaling invariance: 
$x\mapsto x/\lambda$, $t\mapsto 
t/\lambda^{3}$ and $u\mapsto \lambda^{2/(p-1)}u$ with 
$\lambda=const$. For $p=5$, the $L^{2}$ norm of $u$ is invariant 
under this rescaling, and this case is consequently called $L^{2}$-critical. It is shown in \cite{MM} that solutions to the $L^{2}$ 
critical generalised KdV equation can have a blow-up in finite time for 
smooth initial data. The mechanism of the blow-up for initial data 
close to the soliton is discussed in \cite{MMR}. In the present 
paper, we will only study the sub-critical cases, $p\leq 4$.

A convenient way to treat solutions 
varying on length scales of order $1/\epsilon$ for times of order 
$1/\epsilon$ for $0<\epsilon\ll 1$ is to consider the map $x\mapsto 
\epsilon x$, $t\mapsto \epsilon t$. This leads for (\ref{gKdV}) to 
(once more we use the same symbol $u$ for the transformed and the 
original solution)
\begin{equation}
    u_{t}+\epsilon^{2} u_{xxx}+u^{p-1}u_{x}=0. 
    \label{gKdVe}
\end{equation}
The formal limit $\epsilon\to0$ of this equation yields a generalised 
Hopf equation, $ u_{t}+u^{p-1}u_{x}=0$. It is known that such 
equations can have a \emph{hyperbolic blow-up}, i.e., shocks for 
general smooth initial data, for instance for data with a single 
hump. The generic \emph{break-up} of such solutions at a point 
$(x_{c},t_{c},u_{c})$, see for instance 
the discussion in \cite{DGK} and references therein, is characterized 
by the equations
\begin{align}
    a(u_{c})t_{c}+\Phi(u_{c}) & =x_{c},
    \nonumber  \\
    a'(u_{c})t_{c}+\Phi'(u_{c}) &= 0,
    \label{critp}  \\
    a''(u_{c})t_{c}+\Phi''(u_{c}) &= 0,
    \nonumber
\end{align}
where $a(u)=u^{p-1}$. It is known that dispersive regularizations of 
dispersionless equations as (\ref{gKdVe}) in our case will lead to 
\emph{dispersive shock waves} (DSWs), i.e., rapid modulated 
oscillations near the shock of the solution of the corresponding 
dispersionless equation for the same initial data. Dubrovin 
\cite{D06,D08} presented a conjecture that the onset of a DSW is 
\emph{universal} for a large class of dispersive equations and of 
initial data, and that it is given by a special solution of the 
so-called Painlev\'e $P_I^2$ equation, see for instance \cite{kap}. This 
conjecture was numerically shown to apply to the 
generalised KdV equations in \cite{DGK}. In this paper we provide the numerical tools to 
do so for larger classes of initial data than in \cite{DGK}, however 
we leave addressing the universality conjecture in the present context for a future work.

\section{Numerical approach}\label{Sec_Num_appr}
In this section we present the numerical approach to treat 
generalised  KdV equations on the 
compactified real line. We first introduce the compactification which 
allows to study the equation on the interval $[-1,1]$ where we use a 
Chebyshev collocation method. The resulting system of ordinary differential 
equations is then integrated in time with an implicit Runge-Kutta 
scheme. 

\subsection{Compactification}

To numerically treat the KdV equation, we first map the real line via 
the well known map (this is the standard compactification used for 
Minkowski spacetime)
\begin{equation}
    x = c \tan\frac{\pi l}{2}, \quad l\in[-1,1],\quad c=const,
    \label{map}
\end{equation}
to the interval $[-1,1]$. This implies 
\begin{equation}
    \partial_{x}=\frac{2}{\pi c}\cos^{2} \frac{l\pi}{2}\partial_{l}
    \label{derivative}.
\end{equation}

\begin{figure}[htb!]
  \includegraphics[width=0.32\textwidth]{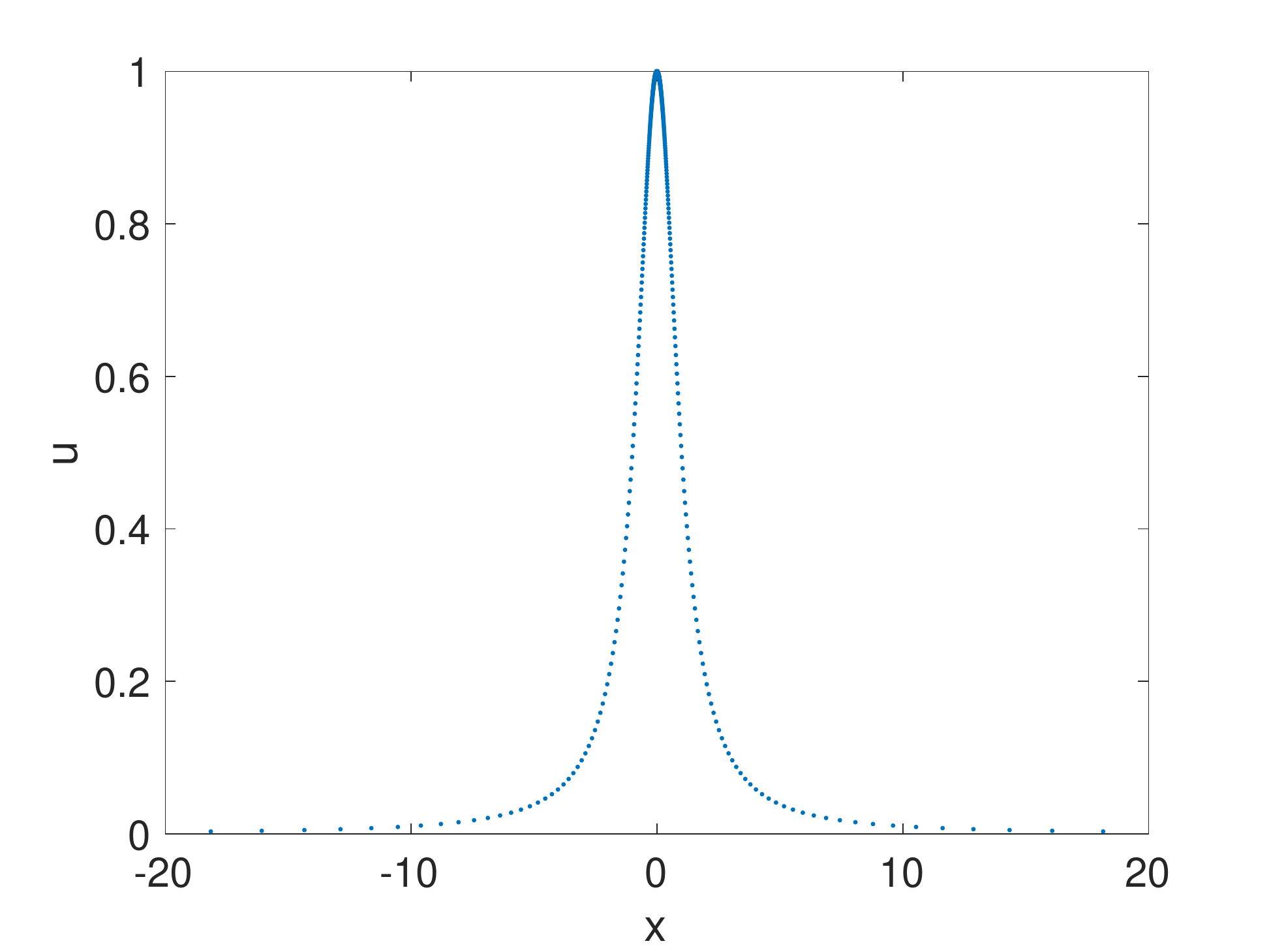}
  \includegraphics[width=0.32\textwidth]{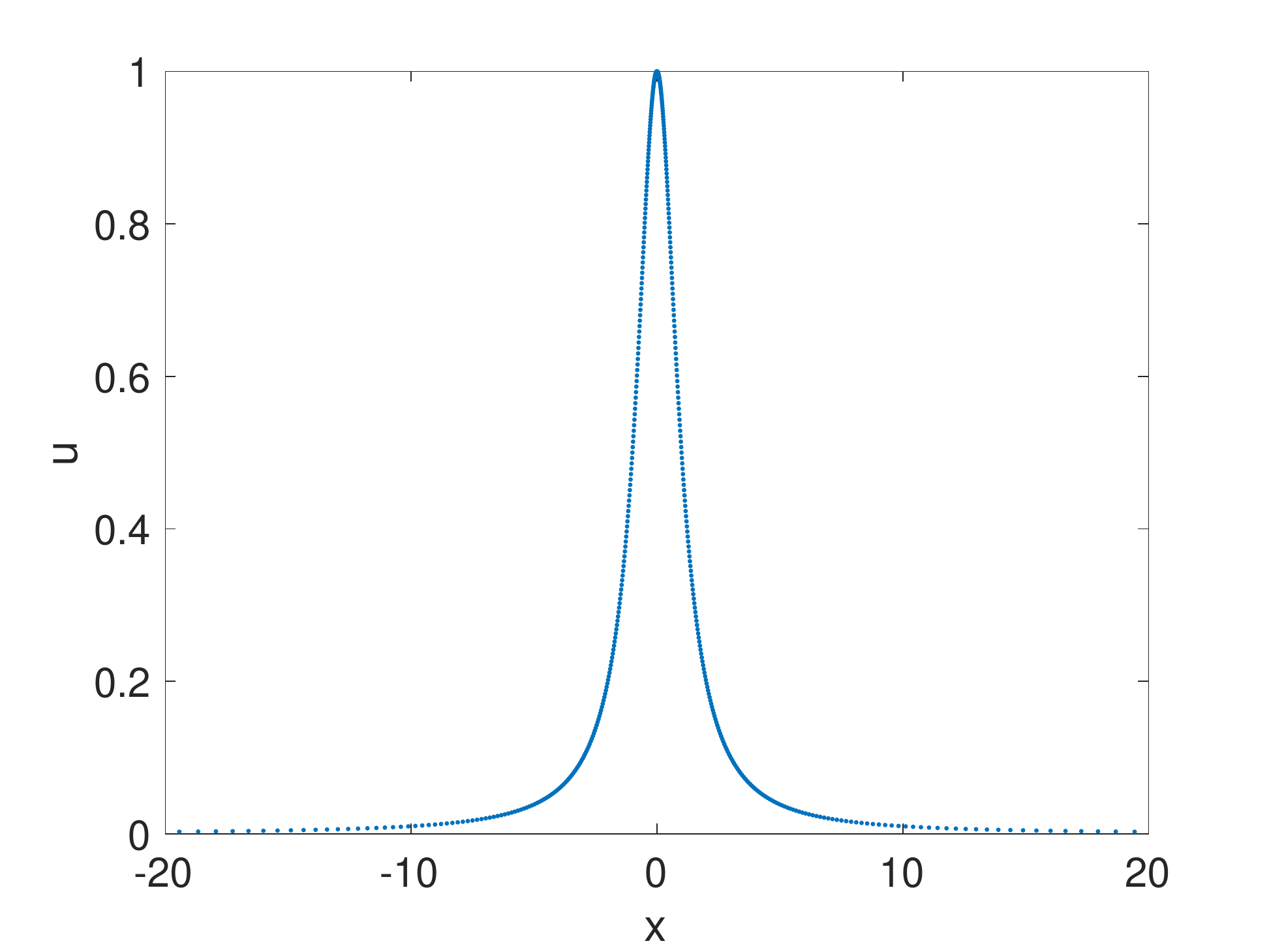}
  \includegraphics[width=0.32\textwidth]{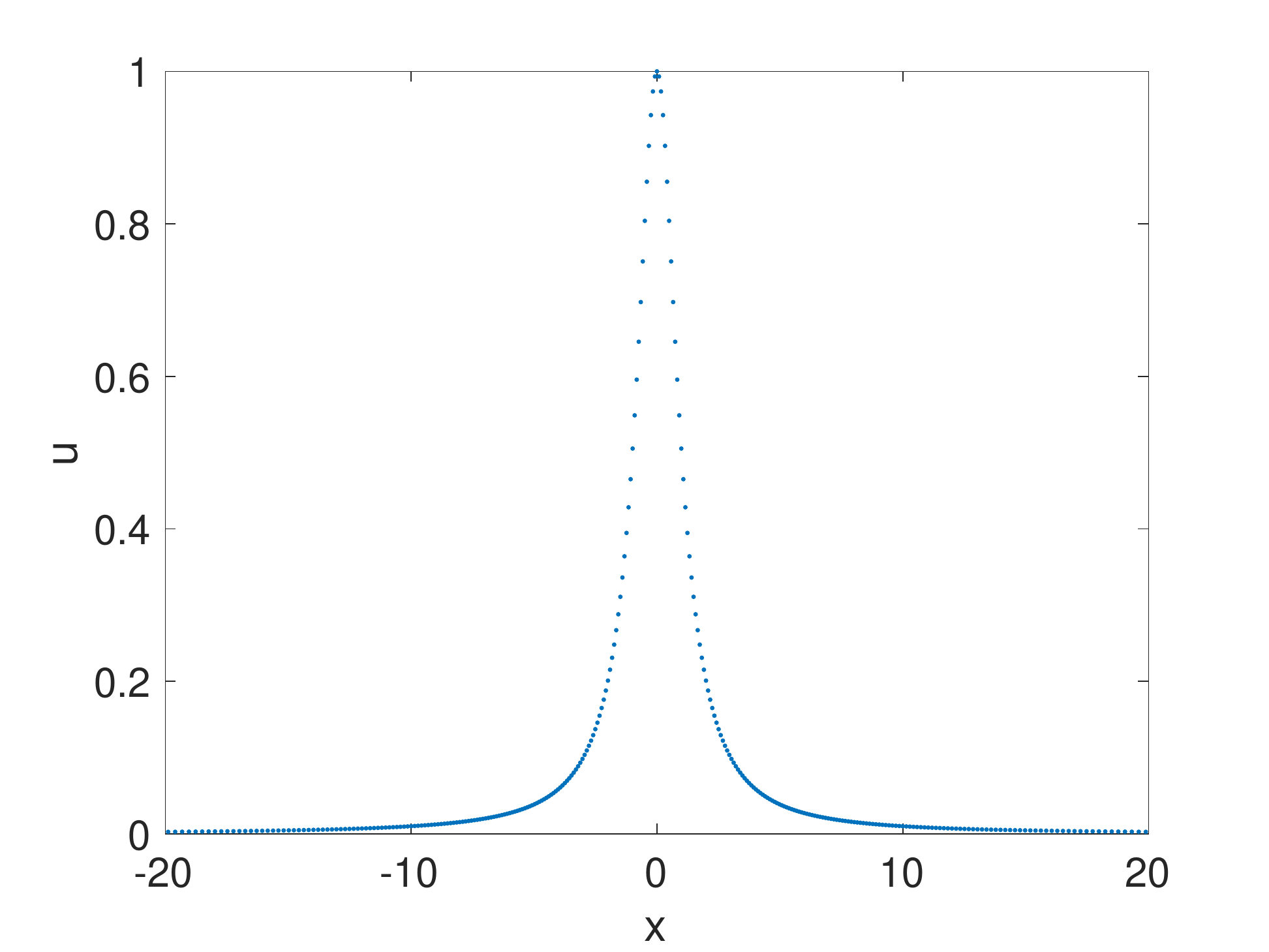}
 \caption{Distribution of the Chebyshev  collocation points (\ref{col}) on a Lorentzian profile under the map (\ref{map}) for $N=800$; on the left for  $c=0.1$, in the middle for $c=1$, on the right for $c=10$.}
  \label{c}
\end{figure}

The role of the constant $c$ is to control the 
numerical resolution in various parts of the real line within certain limits. This is 
illustrated in Fig.~\ref{c} where we have introduced \emph{Chebyshev 
collocation points} for $l\in[-1,1]$, i.e.,
\begin{equation}
    l_{n}=\cos(n\pi/N), \quad n=0,1,\ldots,N, \quad N\in \mathbb{N}
    \label{col}.
\end{equation}
It can be seen in Fig.~\ref{c} that the density of the points near $x=0$ 
is higher for smaller $c$ (numerically the function $\tilde{v}$ is 
studied as a function of $l$, but since in applications the 
dependence on $x$ is more important, we plot its $x$-dependence which 
can be obtained via (\ref{map})). Note, however, that the spectral methods 
we apply in this paper are global. This means that a high resolution 
in part of the studied domain is not necessarily beneficial, only the 
overall resolution on the whole interval is important. Therefore we 
generally apply values of $c$ close to 1. 

\subsection{Boundary conditions}
The map (\ref{map}) transforms the KdV equation (\ref{gKdV}) on 
the real line to an equation on the interval $[-1,1]$, which is 
singular at the points $l=\pm1$. Because of this singular behavior, 
no boundary conditions need to be imposed there. 
However, in practice it is useful to give boundary conditions at 
these points to stabilize the numerical approaches. We apply a
vanishing condition for $u$ at these 
points, and a \emph{clamped boundary condition} for $l=-1$, i.e., the 
three conditions (in an abuse of notation, we denote $u(x,t)$ and 
$u(l,t)$ with the same letter)
\begin{equation}
    u(l,t)\biggr\rvert_{l=1}=0,\quad  u(l,t)\biggr\rvert_{l=-1}=0, 
    \quad  u_l(l,t)\biggr\rvert_{l=-1}=0.
    \label{bound}
\end{equation}

The approach we present here can be generalised to  functions $u$ 
which do not tend to 0 for 
$|x|\to\infty$, but which are bounded there. In this case we write 
\begin{equation}
    u = v + V,\quad V=A\frac{1+l}{2}+B\frac{(1+l)^{2}}{4}
    +C\frac{1-l}{2}
    \label{uv},
\end{equation}
where $A$, $B$, $C$ are constants such that $v(l=\pm1) = 
v_{l}(l=-1)=0$. To treat the clamped boundary condition for $l=-1$, we use as in 
\cite{trefethen} the ansatz
\begin{equation}
    v = (1+l)\tilde{v}. 
    \label{clamp}
\end{equation}

This leads for (\ref{gKdV}) to the equation
\begin{equation}
    \tilde{v}_{t}+\frac{1}{1+l}[(1+l)\tilde{v}+V]_{xxx}
    +\frac{1}{1+l} [(1+l)\tilde{v}+V]^{p}[(1+l)\tilde{v}+V]_{x}=0
    \label{gKdV2}
\end{equation}
which has to be solved for all $t$ with the condition 
$\tilde{v}(\pm1)=0$. 

\subsection{Chebyshev differentiation matrices}\label{chebdisc}
The dependence of $\tilde{v}$ in (\ref{gKdV2}) on $l$ will be treated 
in standard way via Lagrange interpolation of $\tilde{v}$ on Chebyshev collocation 
points (\ref{col}) as discussed in \cite{trefethen}. A derivative of 
$\tilde{v}$ with respect to $l$ is then approximated via the 
derivative of the Lagrange polynomial. This leads to the action of a 
matrix on the vector 
$\tilde{v}=(\tilde{v}(l_{0}),\ldots,\tilde{v}_{N})$ (again we use the 
same symbol for the function $\tilde{v}$ and the vector $\tilde{v}$), the well known 
\emph{Chebyshev differentiation matrices} $D$, see e.g., 
\cite{trefethenweb,chebdiff}. This means with 
(\ref{derivative}) that the 
derivatives $\partial_{x}$ are approximated by 
\begin{equation}
    \partial_{x}\approx \frac{2}{\pi c}\mbox{diag}\left(\cos^{2} \frac{l\pi}{2}
    \right)D
    \label{dxapprox},
\end{equation}
where the diagonal matrix has the components $(\cos^{2} 
\frac{l_{0}\pi}{2},\ldots,\cos^{2} \frac{l_{N}\pi}{2})$.

The Lagrange interpolation of a function on Chebyshev collocation 
points is closely related to an expansion of the function in terms of 
Chebyshev polynomials $T_{n}$, $n=0,1,\ldots$ 
\begin{equation}
    \tilde{v}(l)\approx\sum_{n=0}^{N}v_{n}T_{n}(l),
    \label{coeff}
\end{equation}
see the discussion in \cite{trefethen}. As shown again in \cite{trefethen}, 
the Chebyshev coefficients $v_{n}$, $n=0,1,\ldots,N$ can be computed 
efficiently via a fast cosine transform which is closely related to 
the fast Fourier transform. It is also known that the Chebyshev 
coefficients for a function analytic on $[-1,1]$ decrease 
exponentially, and that the numerical error in approximating a 
function $\tilde{v}$ via (\ref{coeff}) is of the order of the first 
omitted coefficients $v_{n}$ in the Chebyshev series. 

\subsection{Time integration}

After the discretisation in 
space equation (\ref{gKdV2}) becomes 
 an $N+1$-dimensional system of 
ordinary differential equations (ODEs)
in $t$ of the form 
$\tilde{v}_{t}=f(\tilde{v})$ which can be numerically integrated in 
time with standard techniques. The discussion in 
\cite{trefethen} shows that Chebyshev differentiation matrices have a 
conditioning of the order $\mbox{cond}(D^{3})=O(N^{6})$. Thus explicit time 
integration schemes are problematic since stability conditions would 
necessitate prohibitively small time steps. 
Therefore, here we apply an implicit method, and since we are 
interested in capturing rapid oscillations in the expected DSWs, we use a fourth order method. 

Concretely, we apply a fourth  order Runge-Kutta (IRK4) scheme, also 
called the Hammer-Hollingsworth method, a 2-stage 
Gauss scheme. The general formulation of an $s$-stage Runge--Kutta method 
for the initial value problem
$\tilde{v}'=f(\tilde{v},t),\,\,\,\,\tilde{v}(t_0)=\tilde{v}_0$ reads:
\begin{eqnarray}
 \tilde{v}_{n+1} = \tilde{v}_{n} + h      \underset{i=1}{\overset{s}{\sum}} \, 
 b_{i}K_{i}, \\
 K_{i} = f\left(t_{n}+c_ {i}h,\,\tilde{v}_{n}+h  
 \underset{j=1}{\overset{s}{\sum}} \, a_{ij}K_{j}\right),
 \label{K}
\end{eqnarray}
where $b_i,\,a_{ij},\,\,i,j=1,...,s$ are real numbers and
$c_i=   \underset{j=1}{\overset{s}{\sum}} \, a_{ij}$.  
For the IRK4 method used here, one has
$c_{1}=\frac{1}{2}-\frac{\sqrt{3}}{6}$, 
$c_{2}=\frac{1}{2}+\frac{\sqrt{3}}{6}$, $a_{11}=a_{22}=1/4$,
$a_{12}=\frac{1}{4}-\frac{\sqrt{3}}{6}$, 
$a_{21}=\frac{1}{4}+\frac{\sqrt{3}}{6}$ and $b_{1}=b_{2}=1/2$. 

Applying IRK4 to (\ref{gKdV2}) we get the following system,
\begin{align}
    \mathcal{L}K_{1}&=-\frac{1}{1+l}[(1+l)(\tilde{v}+ha_{12}K_{2})+V]_{xxx}\nonumber\\
    &-\frac{1}{1+l} 
    [(1+l)(\tilde{v}+ha_{11}K_{1}+ha_{12}K_{2})+V]^{p}[(1+l)(\tilde{v}+ha_{11}K_{1}+ha_{12}K_{2})+V]_{x},\nonumber\\
        \mathcal{L}K_{2}&=-\frac{1}{1+l}[(1+l)(\tilde{v}+ha_{21}K_{1})+V]_{xxx}\nonumber\\
    &-\frac{1}{1+l} 
    [(1+l)(\tilde{v}+ha_{21}K_{1}+ha_{22}K_{2})+V]^{p}[(1+l)(\tilde{v}+ha_{21}K_{1}+ha_{22}K_{2})+V]_{x},
    \label{K1}
\end{align}
where $$\mathcal{L} = \hat{1}+ha_{11}\frac{1}{1+l}\partial_{xxx}(1+l).$$
Recall that $\partial_{x}$ is approximated in (\ref{K1}) via 
(\ref{dxapprox}), and the third derivative with respect to $x$ is 
approximated as the cubic power of this.

The system (\ref{K1}) will be solved iteratively with a simplified 
Newton iteration. This means that in each step of the iteration the 
new $K_{1}$ and $K_{2}$ are obtained by inverting the operator 
$\mathcal{L}$ only instead of the full Jacobian. The vanishing 
boundary conditions for $\tilde{v}$ and thus for $K_{1}$, $K_{2}$ are 
imposed as in \cite{trefethen}: the equations 
$\mathcal{L}K_{i}=F_{i}$, $i=1,2$, are solved by considering only the 
components $1,\ldots,N-1$ of the $K_{i}$; this means that we consider the reduced 
equations $$
\sum_{m=1}^{N-1}\mathcal{L}_{nm}K_{i,m}=F_{i,n},\quad 
n=0,\ldots,N-1,\quad i=1,2,$$
and solve for $K_{i,1},\ldots,K_{i,N-1}$ only since the values 
$K_{i,0}=K_{i,N}=0$ are imposed. 

The resolution in time can be controlled via the conserved quantities 
of the generalised KdV equation. We consider in general the energy 
for functions vanishing at infinity. For functions which are just 
bounded at infinity, we consider a linear combination of the energy 
and the $L^{2}$ norm of the solution, 
\begin{equation}
    \tilde{E}=\int_{\mathbb{R}}^{} \left( 
    \frac{u^{p+1}-\lambda u^{2}}{p(p+1)}-\frac{1}{2} u_{x}^{2}\right) 
    dx
    \label{Etilde},
\end{equation}
where the constant $\lambda$ is chosen such that the integrand is 
bounded at infinity. The conserved quantities will in actual 
computations depend on time due to numerical errors. But as discussed 
for instance  in \cite{KLE08}, the conservation of these quantities 
controls the numerical error which is in general overestimated by up 
to two 
orders of magnitude.

\section{Examples}\label{Sec_Examples}
In this section we study two types of examples which illustrate the 
potential of the presented numerical approach. First we consider 
initial data in the form of a mollified (smoothed out) step function. Then we study 
initial data not satisfying the Faddeev condition 
$\int^{+\infty}_{-\infty} (1 +|x|)|u(x)|dx <\infty$, see 
\cite{ZAK91}, but nevertheless vanishing for 
$|x|\to\infty$. These examples are studied for two types of 
nonlinearity, for $p=2$, the completely integrable KdV equation, and 
for $p=4$, a non-integrable generalised KdV equation discussed for 
instance in \cite{MM4}. The latter is 
still sub-critical which means that there is no blow-up for 
sufficiently regular initial data.

\subsection{Mollified step initial data}
\begin{figure}[htb!]
  \includegraphics[width=0.49\textwidth]{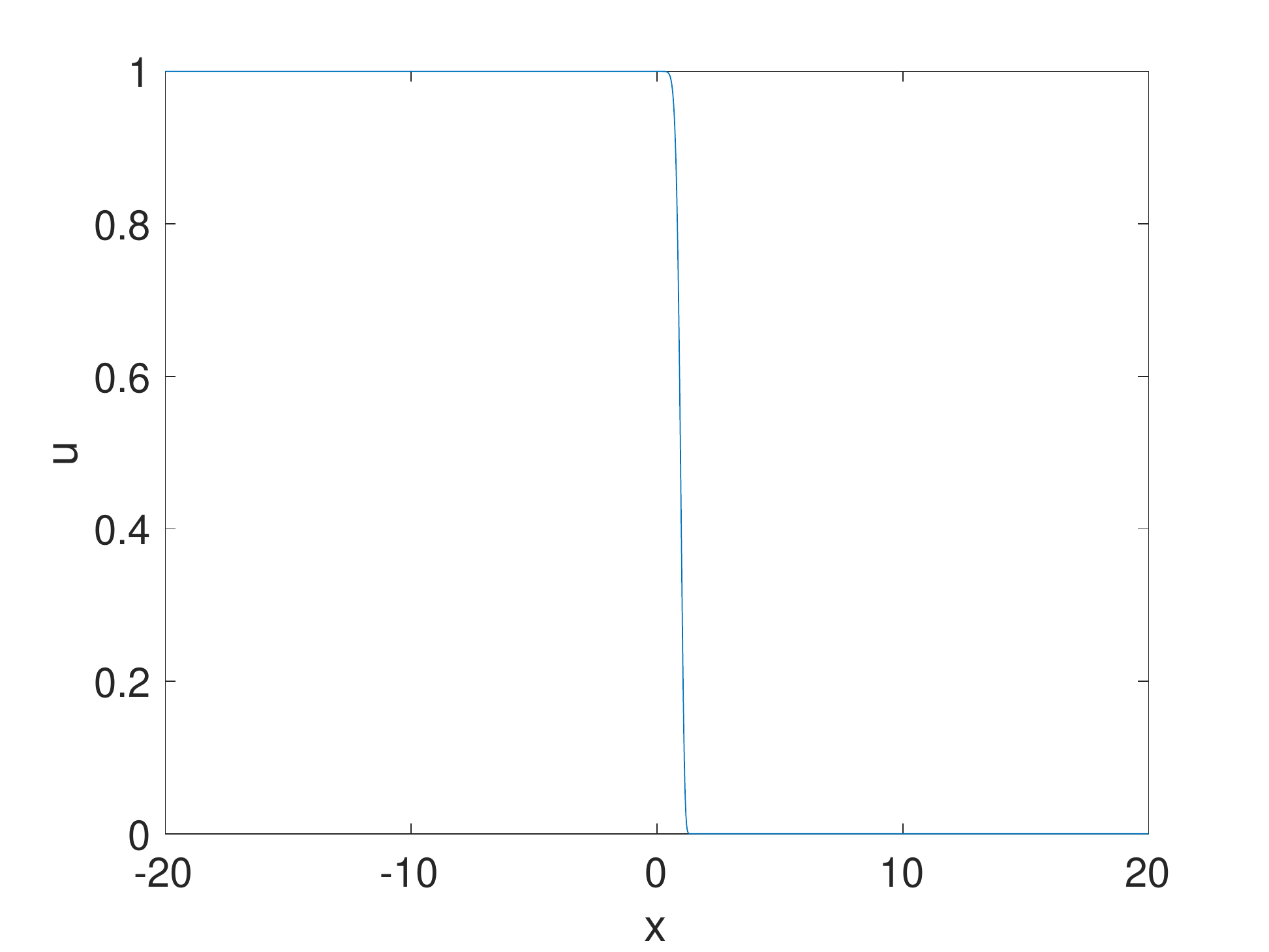}
  \includegraphics[width=0.49\textwidth]{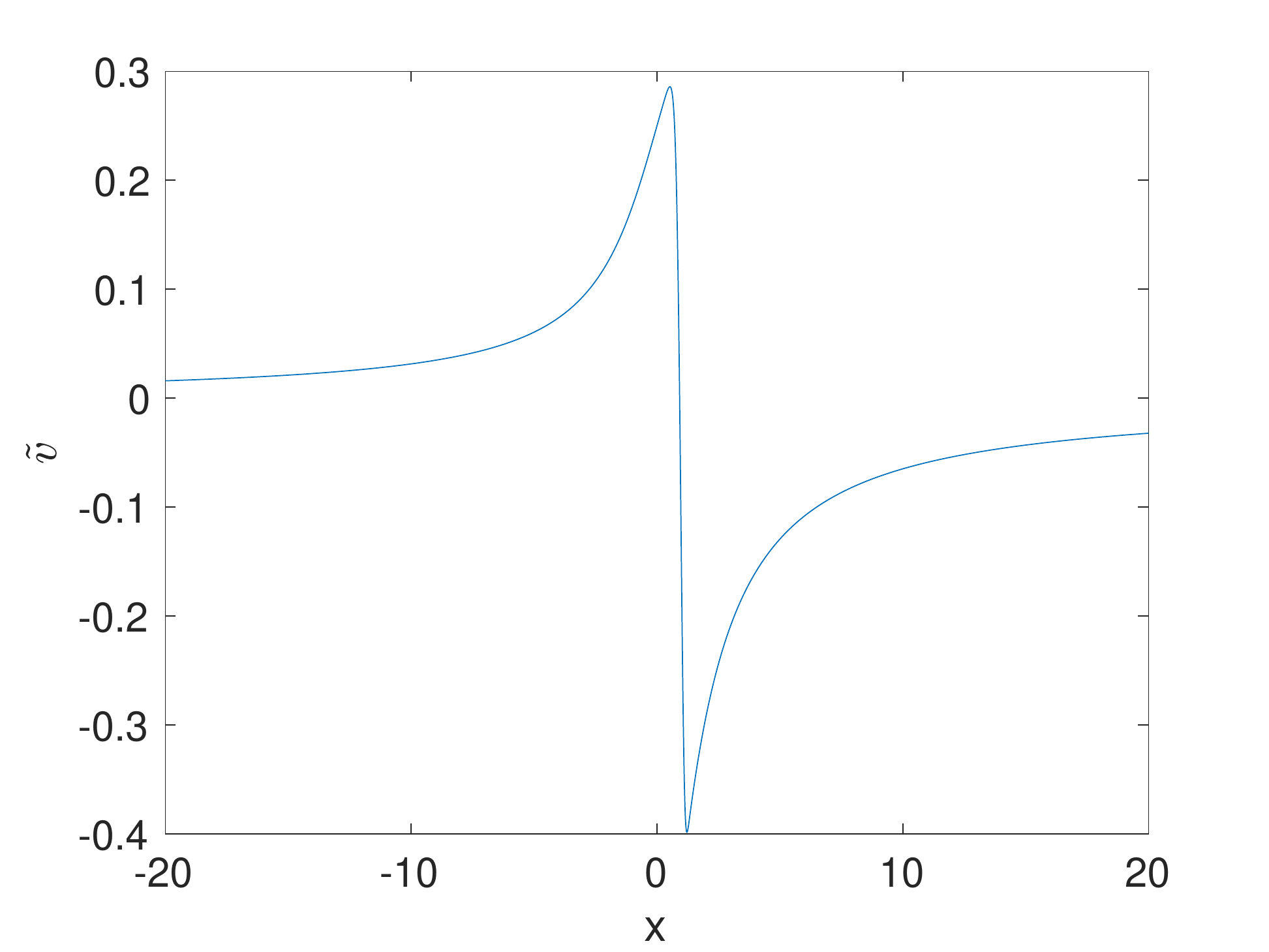}
 \caption{Initial data (\ref{step}) for $n=4$ on the left, and the 
 corresponding  $\tilde{v} = \tilde{v}\left(l(x)\right)$ the right.}
 \label{stepinitial}
\end{figure}

We consider initial data of the form
\begin{equation}
    u(x,0) = 
    \begin{cases}
        1 & x<0 \\
        \exp(-x^{2n}) & x\geq0
    \end{cases}
    \label{step},
\end{equation}
where $n\in\mathbb{N}$. In Fig.~\ref{stepinitial} we show these data 
for $n=4$ on the left, and the corresponding function $\tilde{v}$ on 
the right (one has $A=-B=C=1$). 
\begin{remark}\label{rem}
The function (\ref{step}) is analytic everywhere except at zero, 
where it is $C^{2n-1}$, thus the convergence rate of a spectral 
method is expected to be algebraic. Nevertheless, in practice this 
does not have a detrimental effect on our approach and, as we can see 
below in figure \ref{coefficients_step} on the right,  the behavior 
of Chebyshev coefficients is, due to the finite precision used, virtually the same as in the analytic case.
\end{remark}

\begin{figure}[htb!]
  \includegraphics[width=0.7\textwidth]{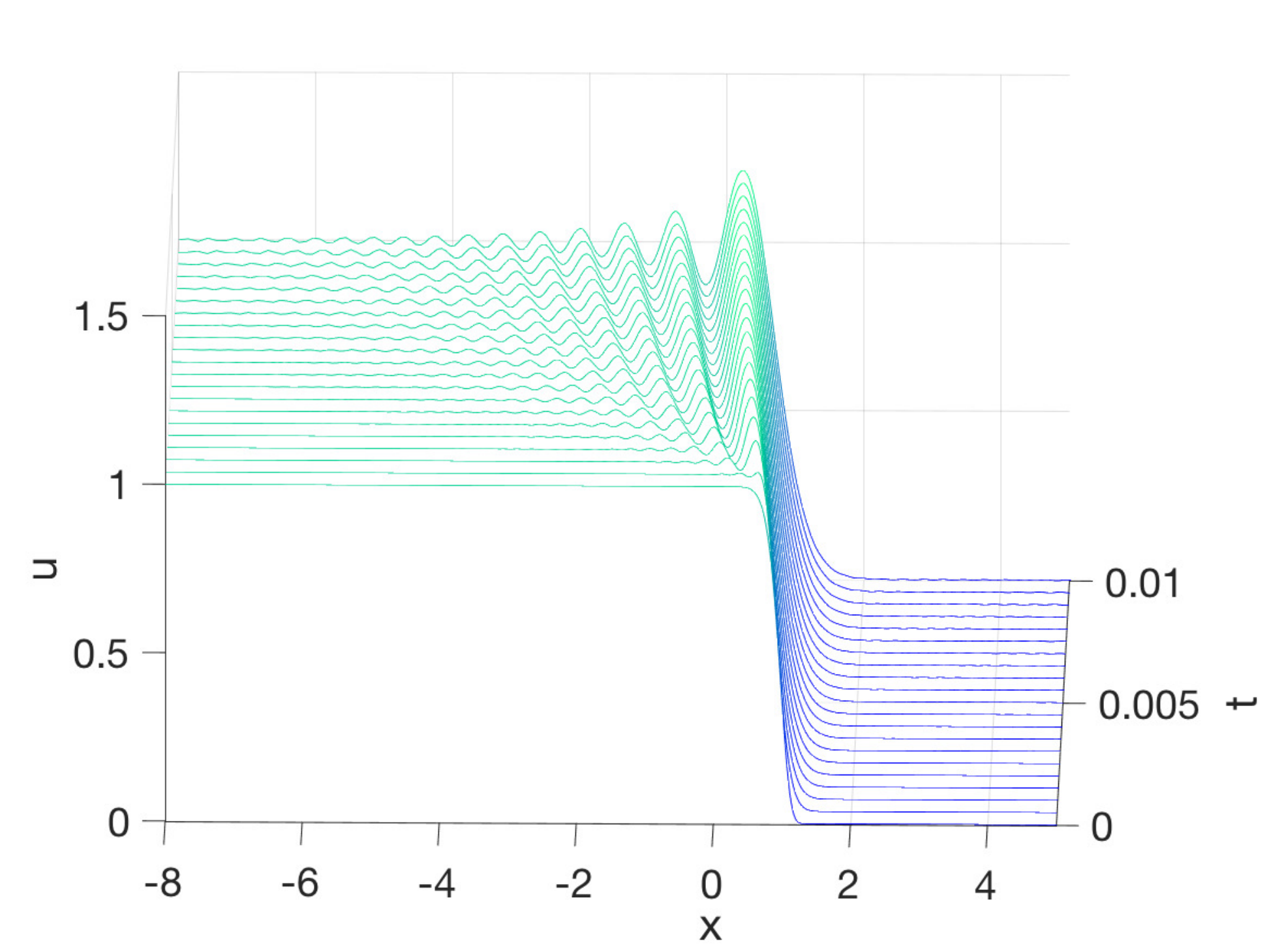}
 \caption{Solution to the KdV equation (\ref{gKdV}) with $p=2$ for 
 the initial data (\ref{step}) for $n=4$ in dependence of time.}
 \label{stepwater}
\end{figure}

For the time evolution of these data, we use $c=2$, $N=600$ and $N_{t}=1000$ 
time steps for $t\in[0,0.01]$. The resulting solution to the KdV 
equation can be seen in Fig.~\ref{stepwater}. The formation of a 
dispersive shock wave is clearly visible. 

\begin{figure}[htb!]
  \includegraphics[width=0.49\textwidth]{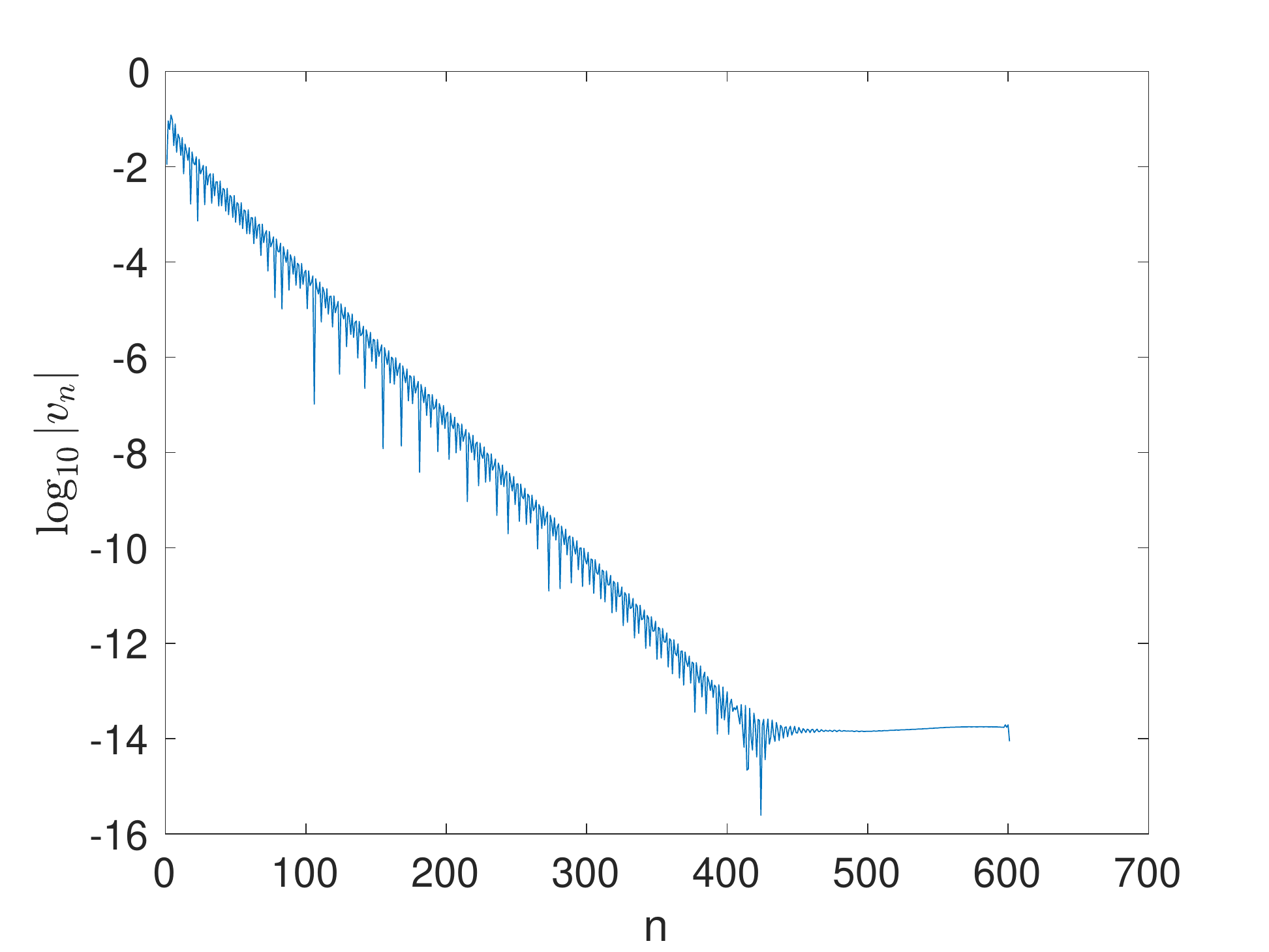}
  \includegraphics[width=0.49\textwidth]{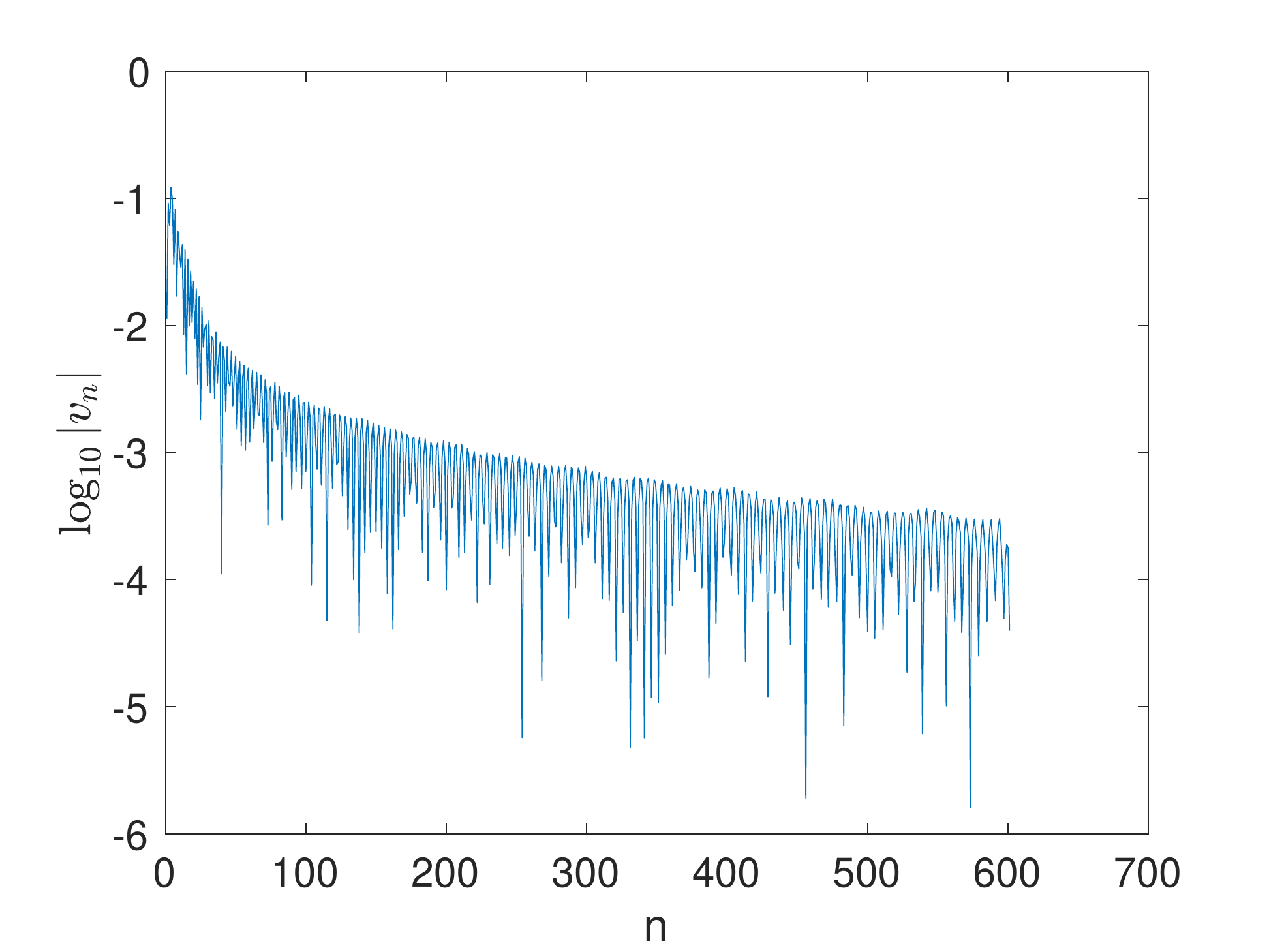}
 \caption{The Chebyshev coefficients (\ref{coeff}), on the left for 
 the mollified step initial data (\ref{step}) 
 on the 
 right for the solution shown in Fig.~\ref{stepwater} for $t=0.01$}
 \label{stepcheb}
\end{figure}

The slow decrease of the amplitude of the oscillations towards 
infinity, similar to the behaviour of the Airy function, is challenging 
for any numerical method. The Chebyshev coefficients $v_{n}$ 
(\ref{coeff})
of the solution  are shown in Fig.~\ref{stepcheb}, on the left for 
$t=0$, on the right for $t=0.01$. They decrease exponentially to the order of the 
rounding error for the initial data indicating that an analytic 
(within numerical precision, see remark \ref{rem})
function is numerically well resolved. The algebraic decay of the 
spectral coefficients for $t=0.01$ indicates an oscillatory 
singularity at infinity as for the Airy function. The spatial 
resolution is thus of the order of $10^{-4}$. The relative 
conservation of the modified energy (\ref{Etilde}) is during the whole computation 
of the order of $10^{-4}$. This means the solution is obtained to 
better than plotting accuracy. 

Note that the DSW is not the same as in the case of an exact step, the 
classical Gurevitch-Pitaevski problem \cite{GP}. But one can for 
instance verify that this is the correct solution by considering a 
finite step smoothed out at both sides, 
\begin{equation}
    u(x,0) = 
    \begin{cases}
        1 & x_{0}<x<0 \\
        \exp(-x^{2n}) & x\geq0\\
        \exp(-(x-x_{0})^{2n}) & x\leq x_{0}
    \end{cases}
    \label{step2}
\end{equation}
which can be conveniently treated with Fourier methods as in 
\cite{KLE08} to which the reader is referred for details and 
references. We use $N=2^{12}$ Fourier modes for $x\in 10[-\pi,\pi]$ and $N_{t}=1000$ time 
steps for a fourth order exponential time differencing method. In 
Fig.~\ref{figstep2} we show on the left the solution of 
Fig.~\ref{stepwater} for $t=0.01$, and for the initial data 
(\ref{step2}) with  $n=4$ and $x_{0}=-5\pi $ at the same 
time.
\begin{figure}[htb!]\label{coefficients_step}
  \includegraphics[width=0.49\textwidth]{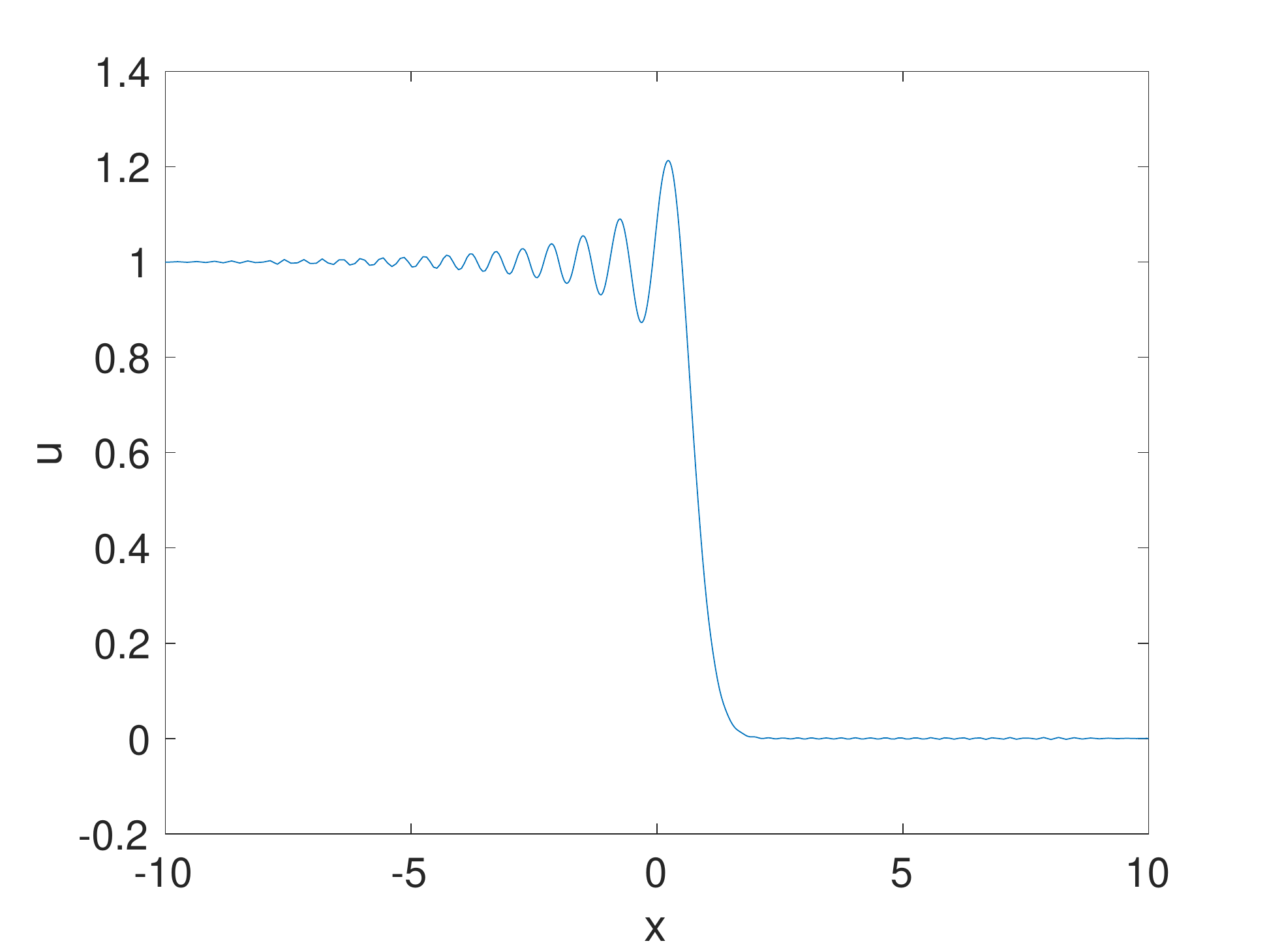}
  \includegraphics[width=0.49\textwidth]{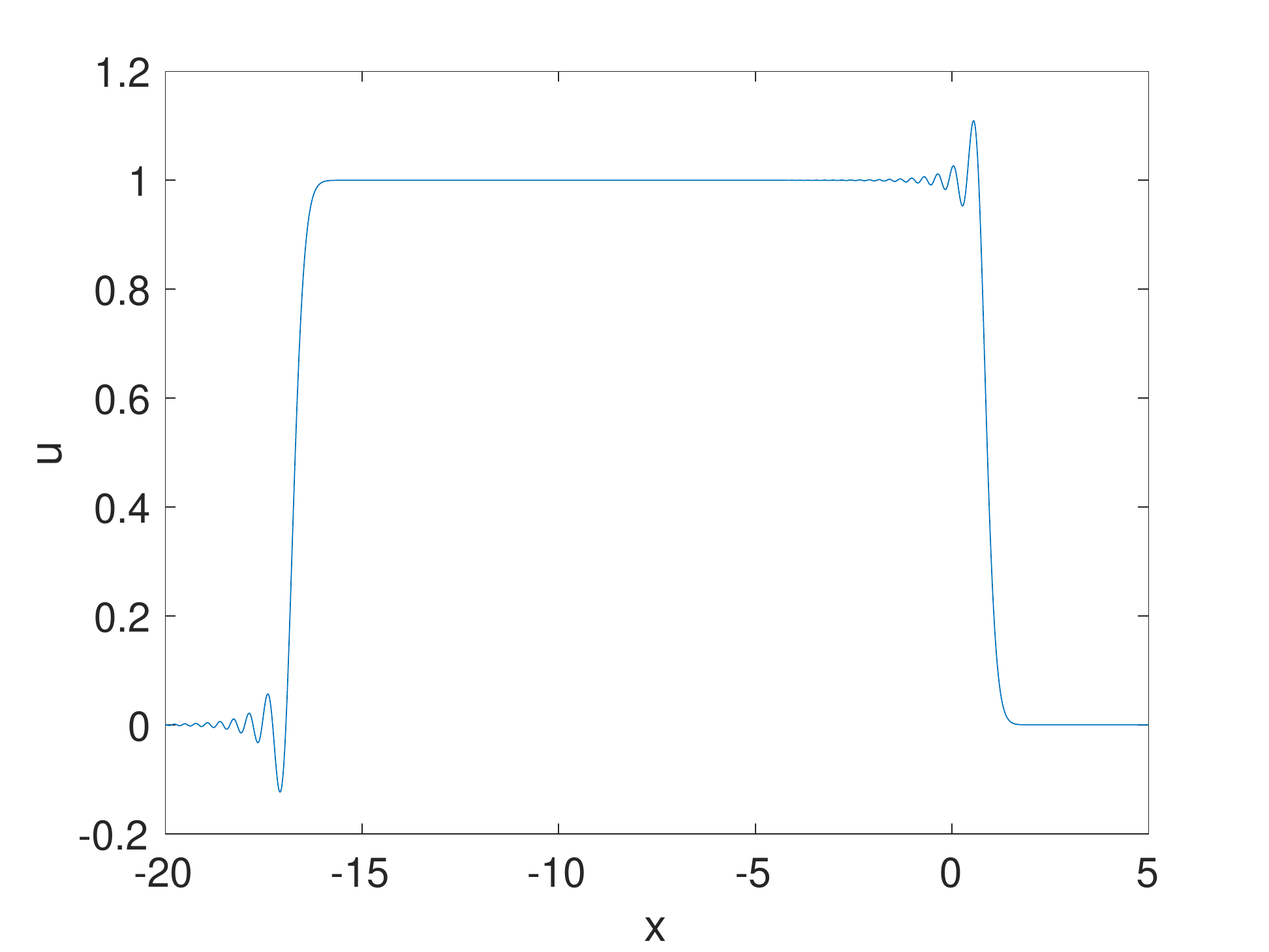}
 \caption{Solution to the KdV equation (\ref{gKdV}) with $p=2$ for 
 the initial data (\ref{step}) on the left, and for the 
 initial data (\ref{step2}) on the right, both for $n=4$ and $t=0.01$.}
 \label{figstep2}
\end{figure}

\begin{figure}[htb!]
  \includegraphics[width=0.7\textwidth]{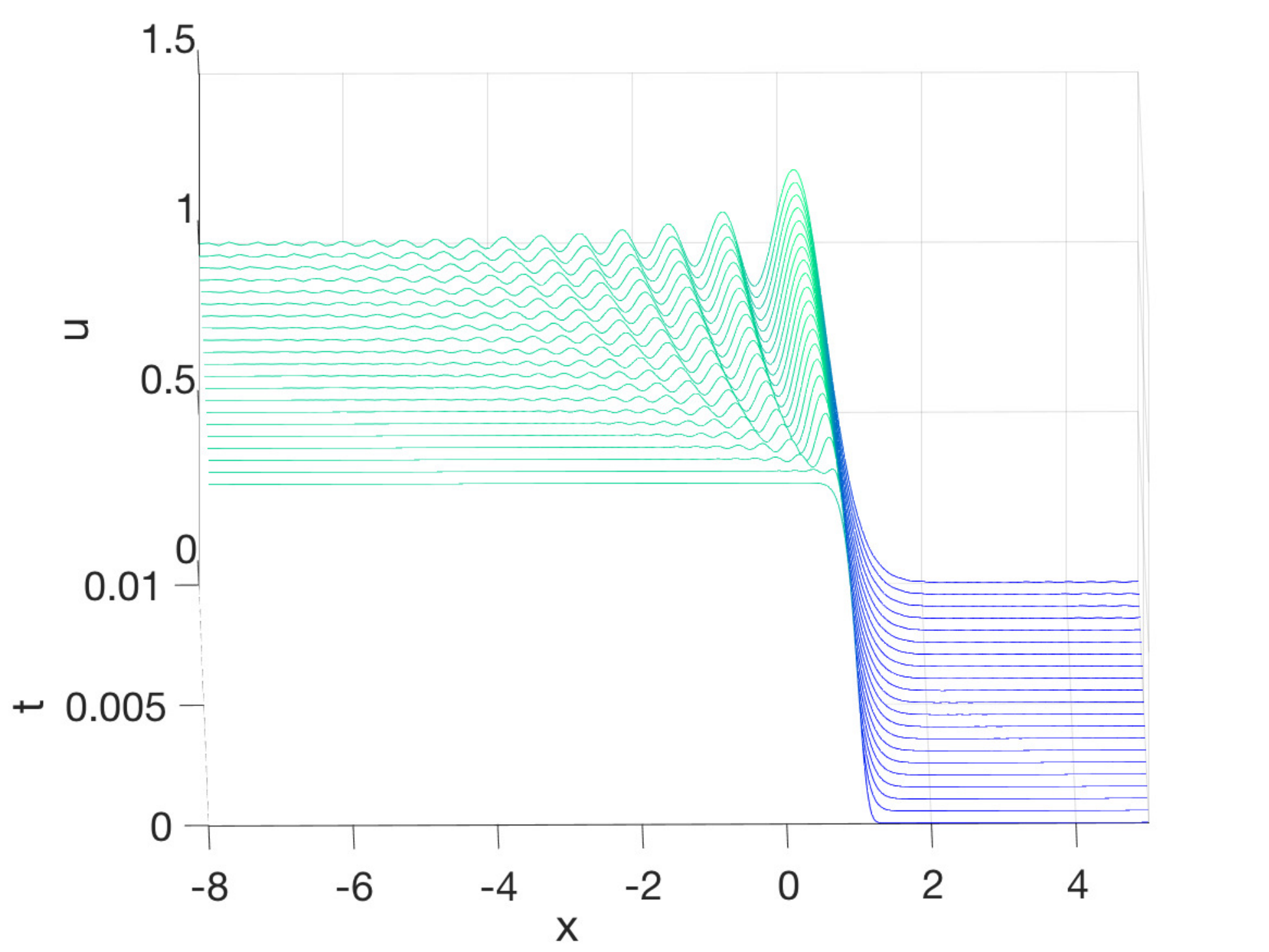}
 \caption{Solution to the generalised KdV equation (\ref{gKdV}) with 
 $p=4$ for 
 the initial data (\ref{step}) for $n=4$ in dependence of time.}
 \label{stepp4water}
\end{figure}

The solution to the generalised KdV equation with $p=4$ for the same 
initial data as in Fig.~\ref{stepinitial} can be seen in 
Fig.~\ref{stepp4water}. We have used the same numerical parameters as 
for the case $p=2$, and we obtain the same numerical resolution. 
The form of the DSW is very similar to one for the 
standard KdV equation. 

\subsection{Slowly decaying initial data}
In this subsection we consider initial data not satisfying the 
Faddeev condition, and we are interested in the long time behavior of 
the corresponding KdV solutions which is done by introducing a small 
parameter $\epsilon$ in (\ref{gKdVe}). Concretely, we study initial data of the form 
\begin{equation}
    u(x,0) = \frac{1}{(1+x^{2})^{a}},\quad a=\frac{1}{2},1
    \label{nfaddeev}.
\end{equation}

\begin{figure}[htb!]
  \includegraphics[width=0.7\textwidth]{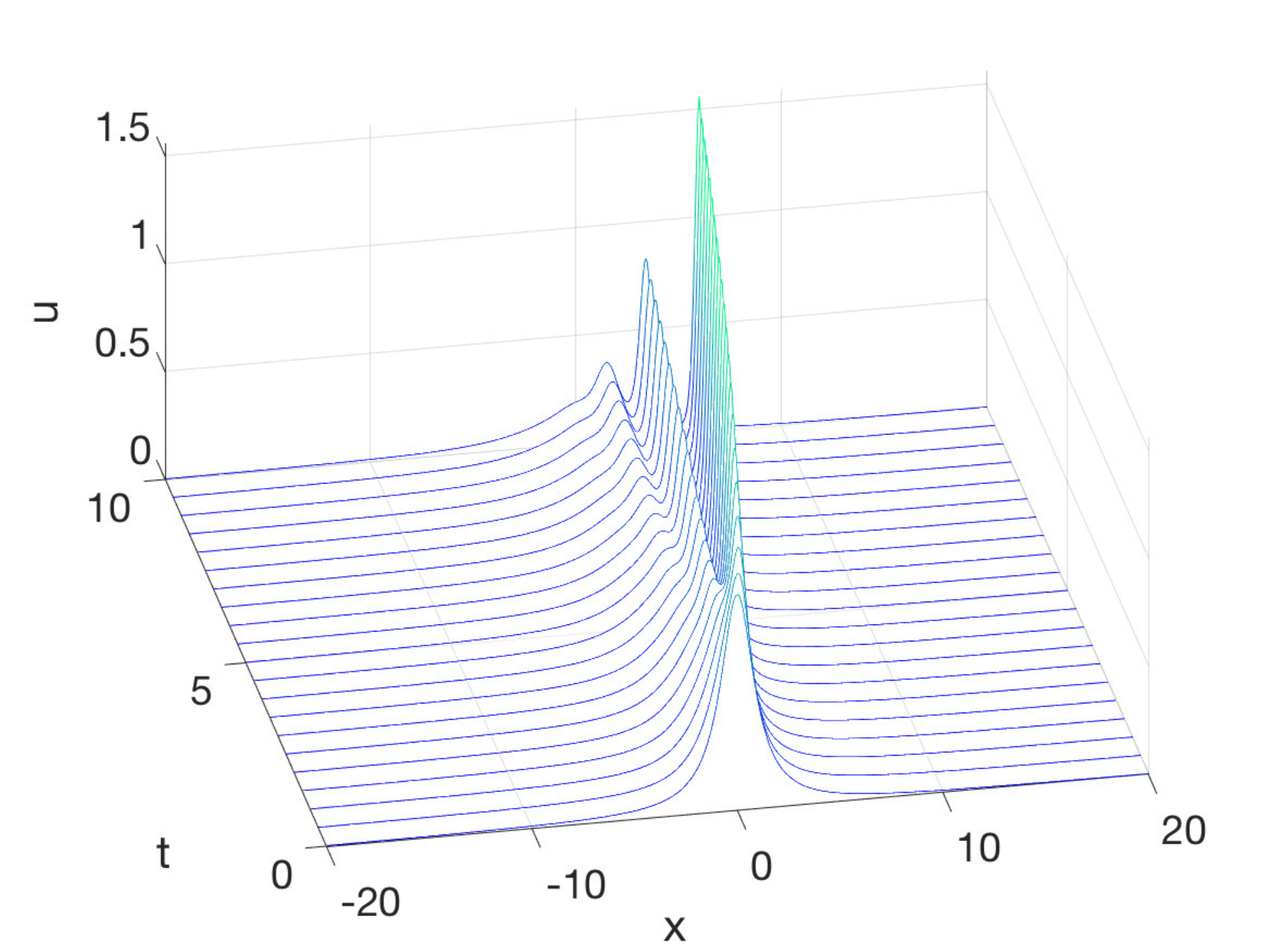}
 \caption{Solution to the KdV equation (\ref{gKdV}) with 
 $p=2$ for 
 the initial data (\ref{nfaddeev}) for $a=1$ in dependence of time.}
 \label{lorentzwater}
\end{figure}

We use $c=2$, $N=800$, and $N_{t}=10^{4}$ time steps for 
$t\in[0,10]$. In Fig.~\ref{lorentzwater} we show the KdV solution 
($p=2$) in (\ref{gKdVe}) for the initial data (\ref{nfaddeev}) for 
$\epsilon=10^{-1}$ and $a=1$ in dependence of time. It can be seen that 
several solitons appear.

\begin{figure}[htb!]
  \includegraphics[width=0.7\textwidth]{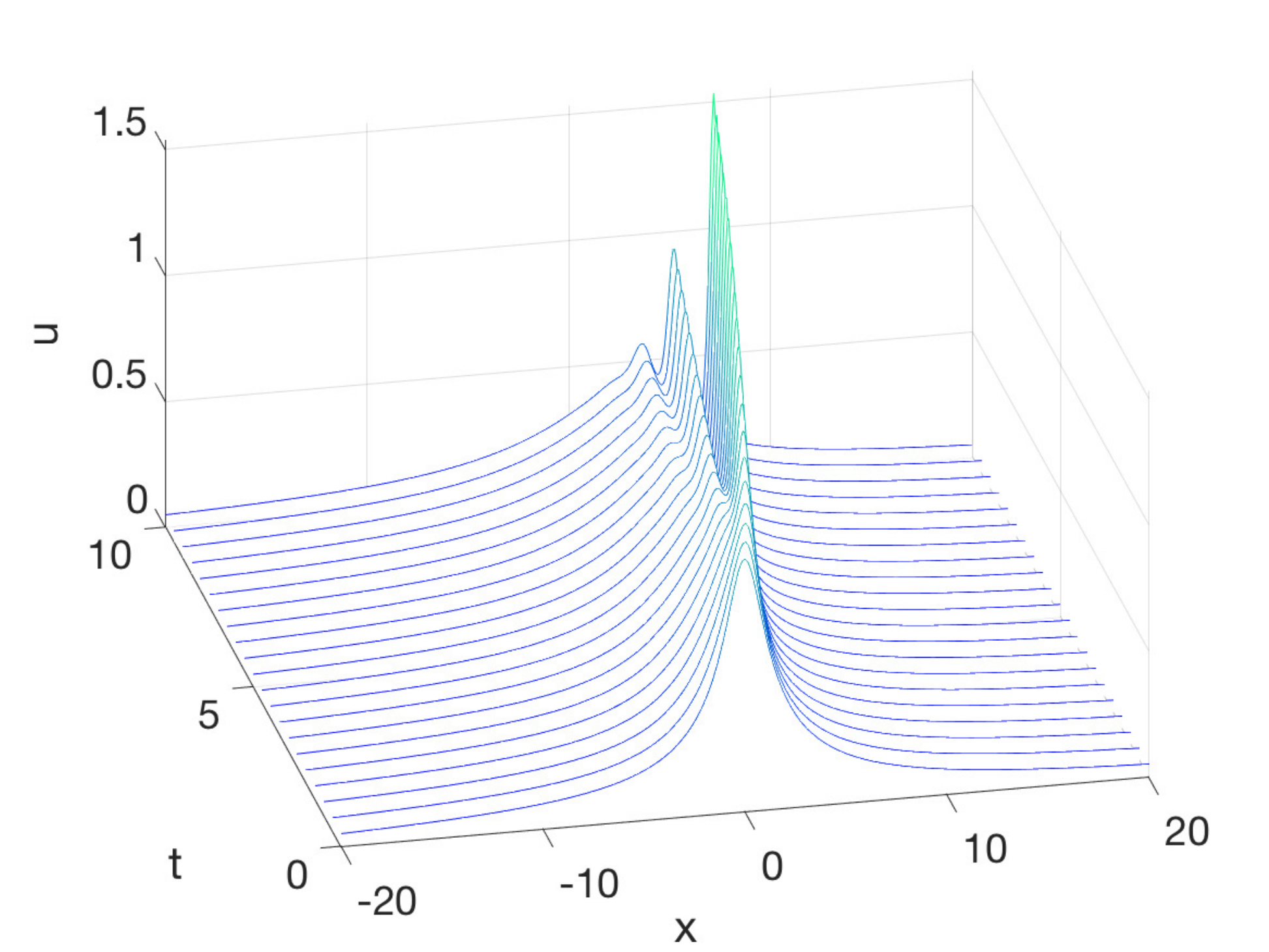}
 \caption{Solution to the KdV equation (\ref{gKdV}) with 
 $p=2$ for 
 the initial data (\ref{nfaddeev}) for $a=1/2$ in dependence of time.}
 \label{rootwater}
\end{figure}

The corresponding KdV solution for the initial data 
(\ref{nfaddeev}) with $a=1/2$  and $\epsilon=10^{-1}$ can be seen in Fig.~\ref{rootwater}.   
Note that in contrast to the case $a=1$, the initial data do not satisfy 
the clamped boundary conditions for $x\to-\infty$, one has 
$A=-B=\pi/2$ and $C=0$ in (\ref{uv}). 

\begin{figure}[htb!]
  \includegraphics[width=0.49\textwidth]{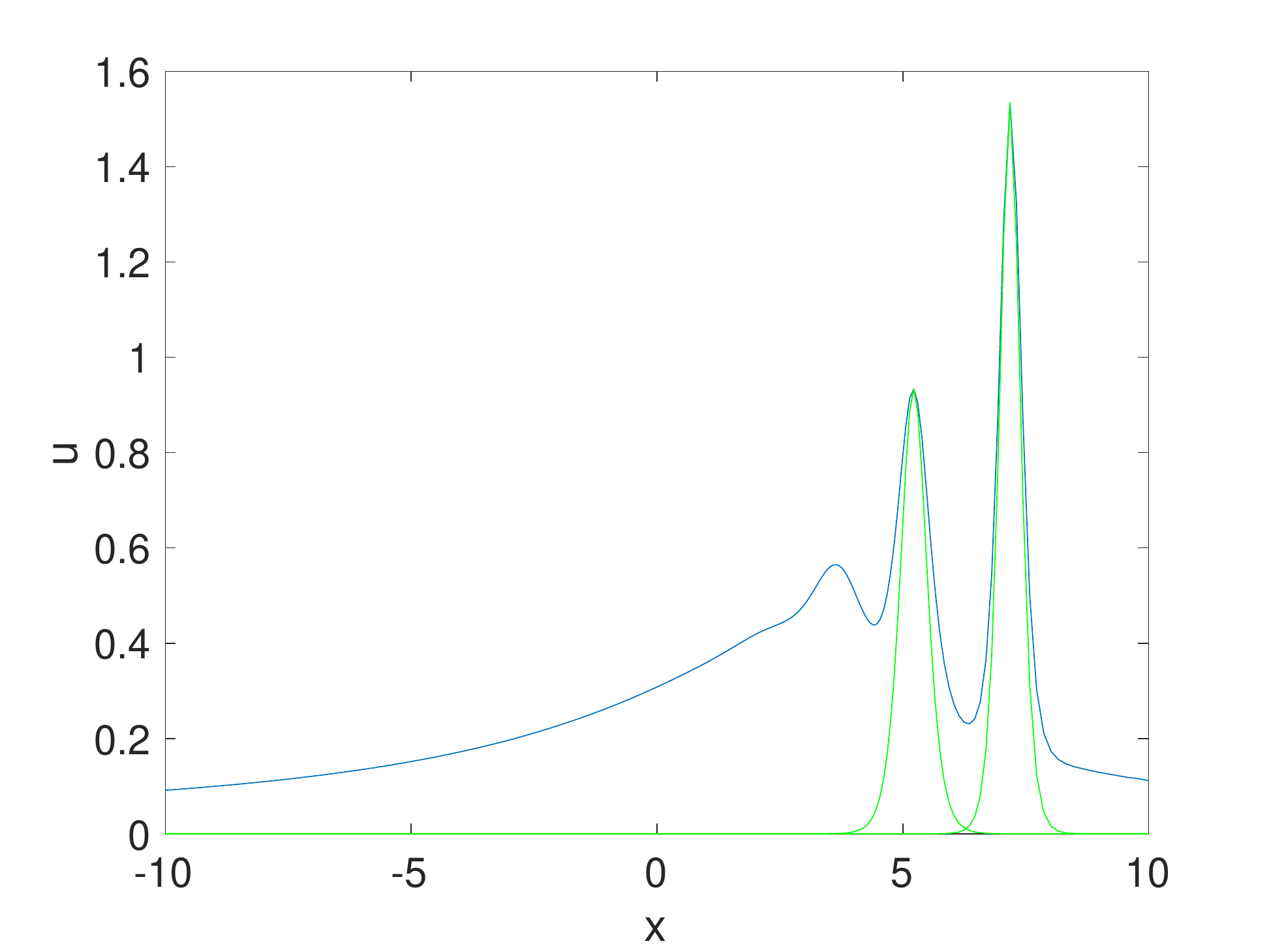}
  \includegraphics[width=0.49\textwidth]{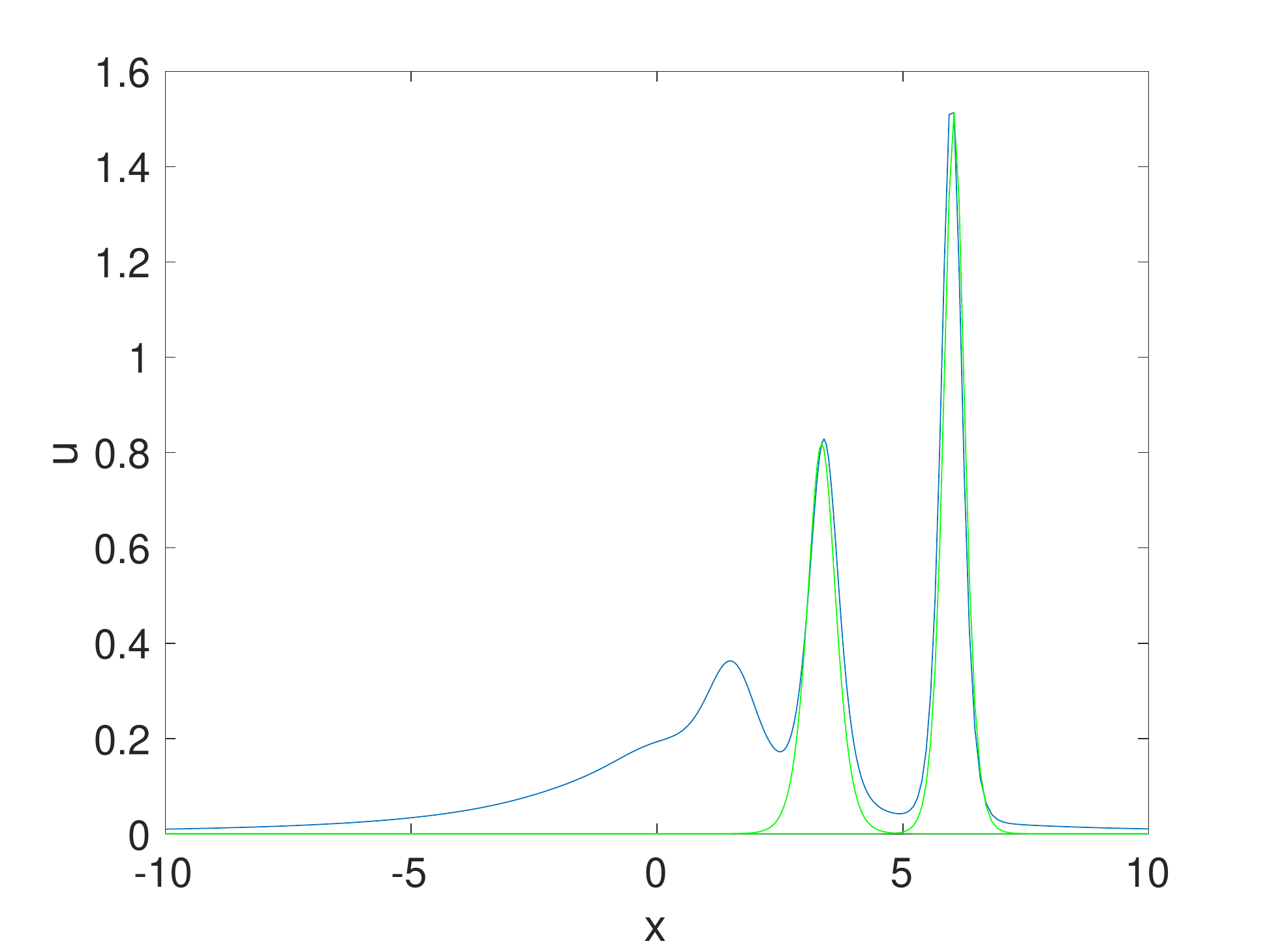}
 \caption{Solution to the KdV equation (\ref{gKdV}) with 
 $p=2$ for 
 the initial data (\ref{nfaddeev}) for $t=10$, on the left for $a=1/2$, 
 on the right for $a=1$; in green fitted solitons (\ref{soliton}).}
 \label{2a}
\end{figure}

The solutions at the final time of Fig.~\ref{lorentzwater} and 
\ref{rootwater} can be seen in Fig.~\ref{2a}, on the left for 
$a=1/2$, on the right for $a=1$. The slower decay towards infinity of 
the initial data with $a=1/2$ can be recognized. But at the time 
$t=10$, one observes the same number of solitons in both cases. The 
peaks in the solutions have been fitted to the solitons 
(\ref{soliton}) which are shown in green in the same figure. It can 
be seen that the solitons are not yet fully separated from the 
background, but that they can be already clearly identified at this 
early stage. 

\begin{figure}[htb!]
  \includegraphics[width=0.49\textwidth]{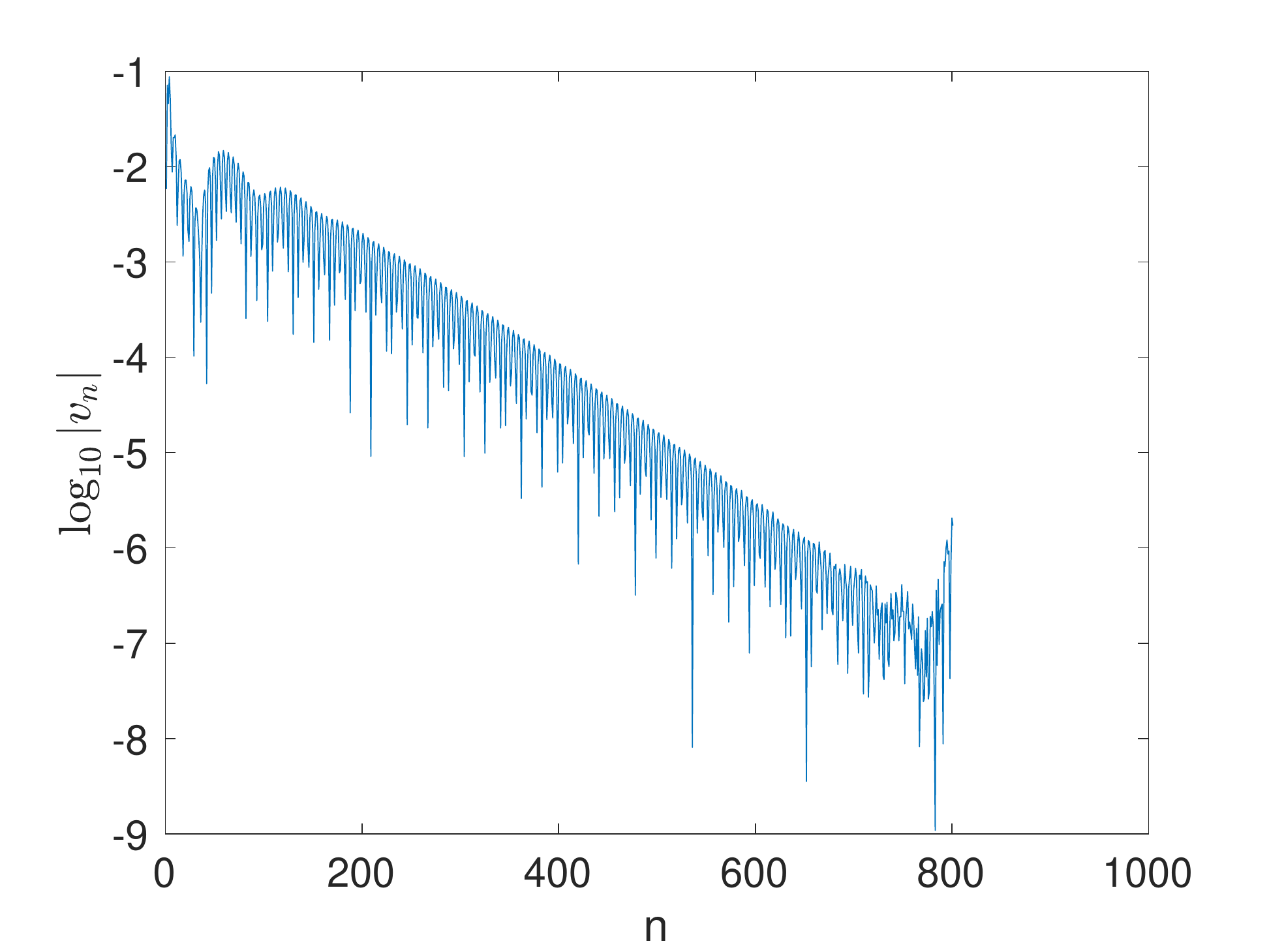}
  \includegraphics[width=0.49\textwidth]{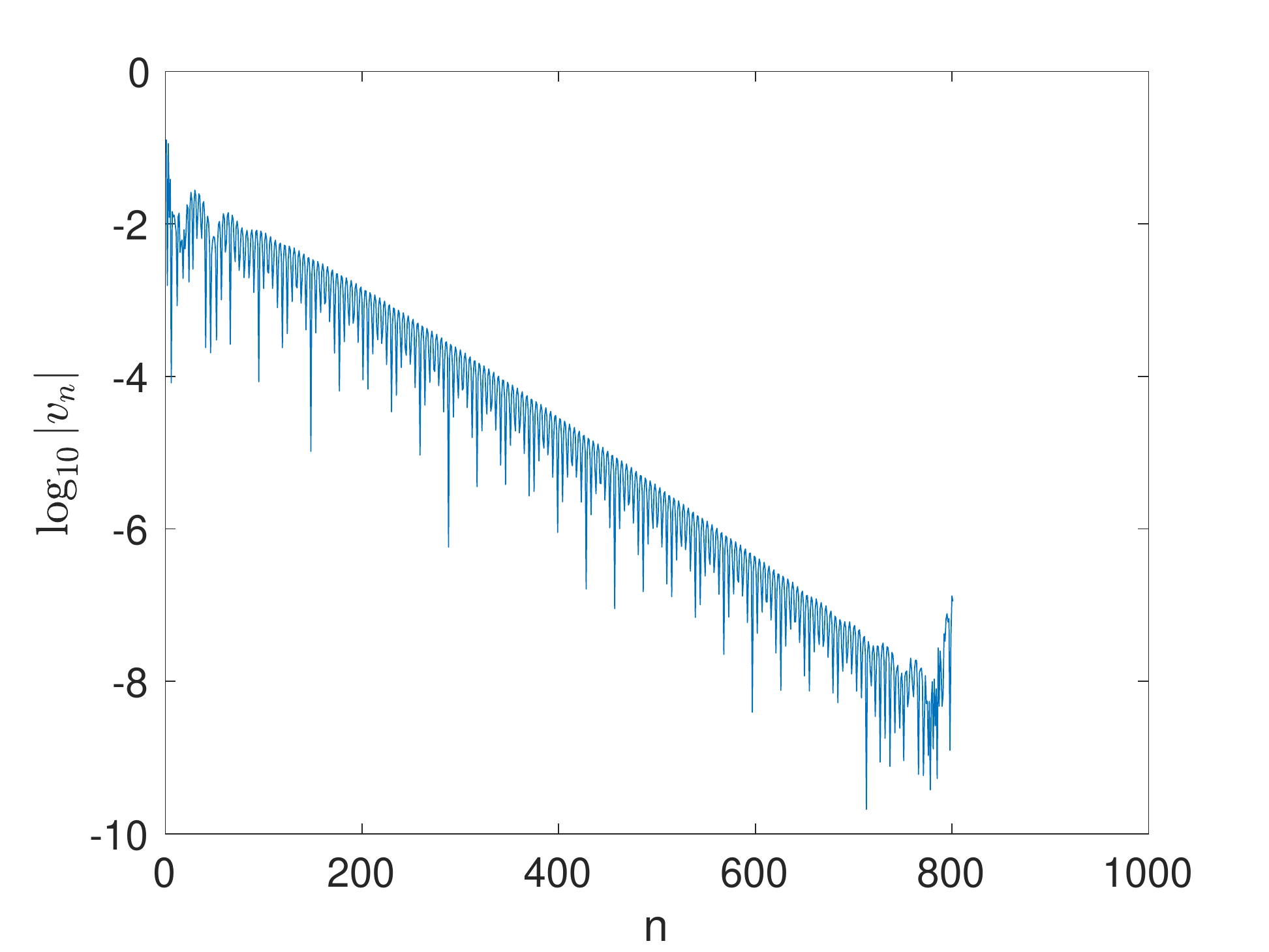}
 \caption{Solution to the KdV equation (\ref{gKdV}) with 
 $p=2$ for 
 the initial data (\ref{nfaddeev}) for $t=10$, on the left for $a=1/2$, 
 on the right for $a=1$.}
 \label{2acheb}
\end{figure}

The relative computed energy is in both cases conserved to the order 
of $10^{-10}$. The Chebyshev coefficients for $t=10$ are shown in 
Fig.~\ref{2acheb}, on the left for $a=1/2$, on the right for $a=1$. 
It can be seen that the coefficients decrease as expected 
exponentially, and that the solutions are well resolved in space as well.   

\begin{figure}[htb!]
  \includegraphics[width=0.49\textwidth]{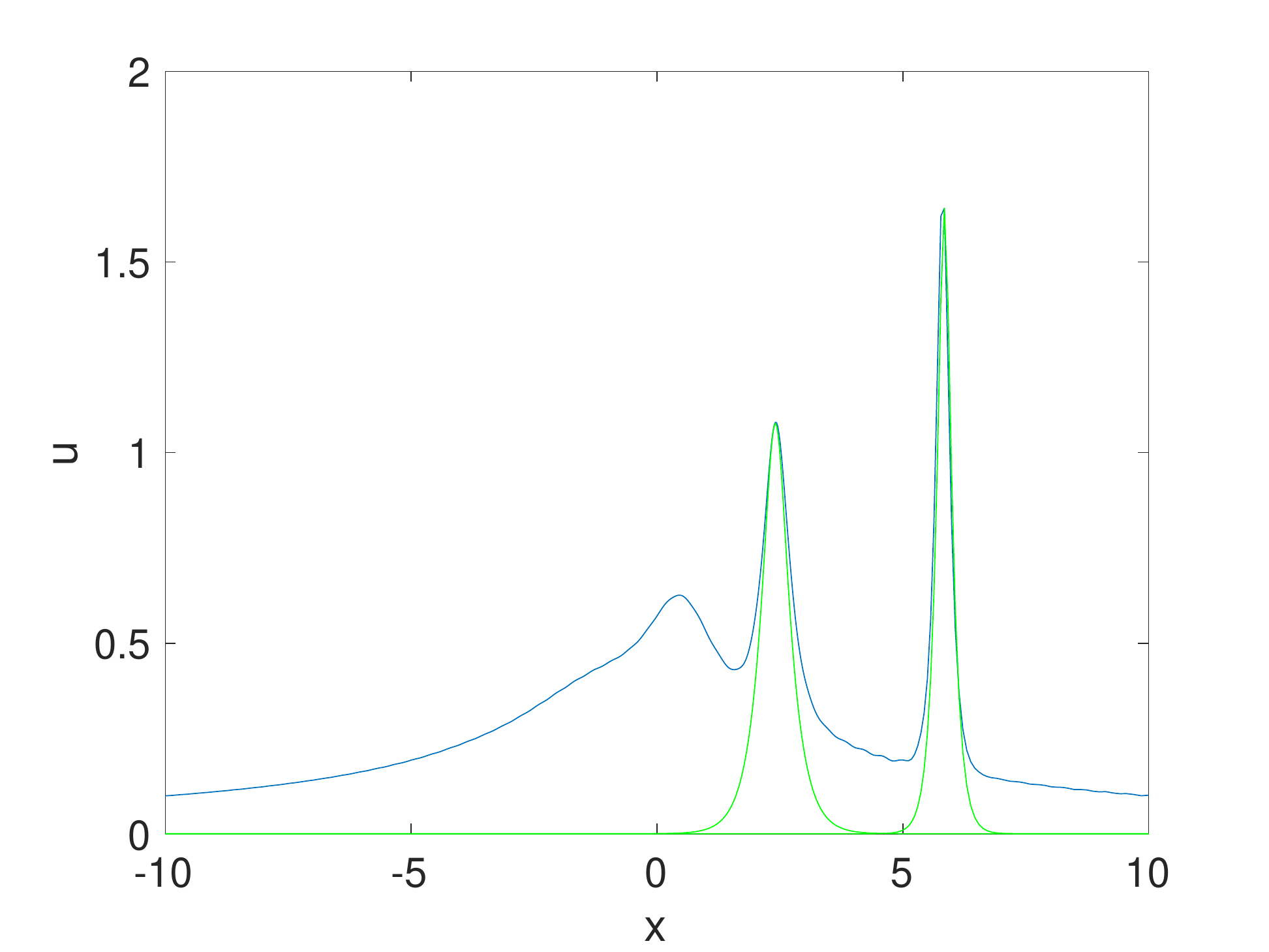}
  \includegraphics[width=0.49\textwidth]{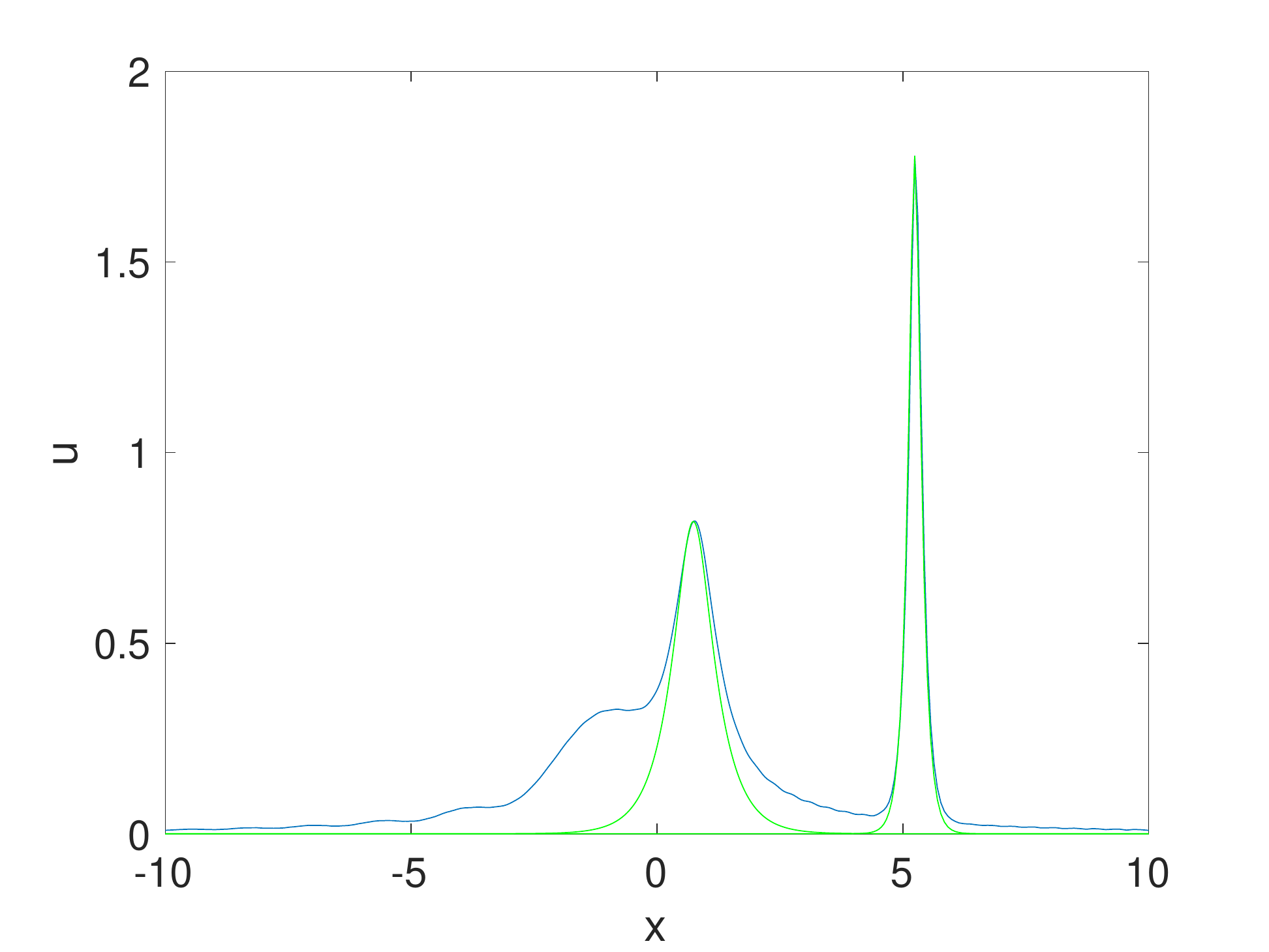}
 \caption{Solution to the KdV equation (\ref{gKdV}) with 
 $p=4$ for 
 the initial data (\ref{nfaddeev}) for $t=10$, on the left for $a=1/2$, 
 on the right for $a=1$.}
 \label{2ap4}
\end{figure}

If the same initial data as in Fig.~\ref{2a} are considered for the 
generalised KdV equation (\ref{gKdVe}) with $p=4$, one obtains for $t=10$ 
the solutions shown in Fig.~\ref{2ap4}. The solitons are here much 
more peaked than in the KdV case of Fig.~\ref{2a} which is also 
illustrated by the fit to the solitons. Consequently the 
same numerical parameters as there lead for the generalised KdV 
solution to a lower resolution: the relative conservation of the 
energy is of the order of $10^{-4}$, and Chebyshev coefficients 
decrease still exponentially, but only to the order of $10^{-6}$ in 
this case.  

\section{Outlook}\label{Sec_Conclusion}
In this paper we have presented  a numerical approach for generalised 
KdV equations on the compactified real line which allows to 
approximate functions which are smooth on $\mathbb{R}\cup\{\infty\}$ 
with spectral accuracy, i.e., with a numerical error decreasing 
exponentially with the number of collocation points. The time 
integration is performed with an implicit fourth order method. One 
direction of further research will be to improve the efficiency of 
the time integration, ideally an explicit approach, for instance 
similar to the approach of \cite{KS18} in the context of the 
Schr\"odinger equation. 

Of special interest is, however, the application of the techniques of 
the present paper to the numerical study of blow-up in the context of 
generalised KdV equations, i.e., for (\ref{gKdV}) with $p\geq 5$. The 
current approach would allow to study the dynamically rescaled 
generalised KdV equations, see for instance \cite{KP13} without the 
problems there related to the use of Fourier methods. In this context 
it would be beneficial to apply a multidomain spectral method as in 
\cite{birem} for Schr\"odinger equations where the compactified real axis is divided into several domains each of which is mapped to the interval 
$[-1,1]$. This allows for a more efficient allocation of numerical 
resolution than via the choice of the parameter $c$ in (\ref{map}), 
in particular a higher resolution near the expected blow-up. To this 
end, matching conditions at the domain boundaries (the solution $u$ 
has to be $C^{2}$ there) have to imposed which is numerically 
problematic. This will be a subject of further research.

\end{document}